\documentclass[letterpaper,12pt]{article}
\usepackage{amsmath}
\usepackage{amssymb}
\usepackage{geometry}
\usepackage[colorlinks=true,linkcolor=red,citecolor=blue]{hyperref}
 \usepackage{amssymb}
 \usepackage{amsmath}
 \usepackage[latin1]{inputenc}
 \usepackage[T1]{fontenc}
\usepackage{times}
\usepackage{epsfig}
 \usepackage{dsfont}
 \usepackage{natbib}
 \usepackage{mathrsfs}
\date{\today}

\newtheorem{theorem}{Theorem}[section]

\newtheorem{corollary}[theorem]{Corollary}

\newtheorem{remark}[theorem]{Remark}

\numberwithin{equation}{section}

\begin{document}

\title{Strong Approximation of Empirical Copula Processes by Gaussian Processes}
\author{Salim BOUZEBDA\footnote{e-mail: salim.bouzebda@upmc.fr~~(Corresponding author)}\hbox{ }
and Tarek ZARI\footnote{e-mail: zaritarek@gmail.com} \\L.S.T.A.,
Universit\'e Pierre et Marie Curie\\ 4 place Jussieu
    75252 Paris Cedex 05, France}
\date{\today}
\maketitle

\begin{abstract}
\noindent
We provide the strong approximation
of normalized empirical copula process by a Gaussian process. In addition we establish a strong approximation
of the smoothed empirical copula processes and a law of iterated logarithm.

%This paper investigates the problem of strong approximation of the
%empirical copula processes for arbitrary dimension, with continuous
%unknown margins. The idea of the proof is based on the results obtained in
%\cite{deheuvels2006}, and the theorem of strong
%approximation for an arbitrary distribution function proved in
%\cite{Csorgo1988}. Using these results, we derive the normality
%for smoothed empirical copula processes and the L.I.L. for empirical copula processes.

 \end{abstract}
 \section{Introduction}
The aim of the present paper is to provide the strong  approximations of the \textit{normalized empirical copula process}
$\{\mathbb{A}_n(\mathbf{u}) : \mathbf{u} \in [0, 1]^d,~~n\geq 1\}$, (see, e.g., (\ref{processuscopule})
below for definition), by a single Gaussian process
$\{\mathscr{K}^{*}_{\mathbb{ C}}(\mathbf{u},n) : \mathbf{u} \in [0, 1]^d,~~n\geq 1\}$.
Thus, we get the strong approximations of $\{\mathbb{A}_n(\mathbf{u}) : \mathbf{u} \in [0, 1]^d,~~n\geq 1\}$ in
terms of Gaussian process in both $\mathbf{ u}$ and $n$. We will
 mainly be concerned with the general case, in which $\{\mathbb{A}_n(\mathbf{u}) : \mathbf{u} \in [0, 1]^d,~~n\geq 1\}$ is generated by
a sample of random vectors with dependent marginals.\\
Consider a continuous random vector $\mathbf{X}=(X_{1},\ldots, X_d )$ with joint cumulative distribution function $\mathbb{ F}(\mathbf{x})=
\mathbb{P}(\mathbf{ X}\leq \mathbf{ x})=
\mathbb{P}(X_{1}\leq x_{1},\ldots,X_{d}\leq x_{d})$, for $\mathbf{ x}\in \mathbb{R}^{d}$, and margins $F_1(\cdot),\ldots, F_{d}(\cdot)$.
Here and elsewhere, for $\mathbf{ x}=(x_{1},\ldots,x_{d})$ and $\mathbf{ y}=(y_{1},\ldots,y_{d})$, we write $\mathbf{ x}\leq \mathbf{ y}$
to denote that $x_{j}\leq y_{j}$, for $j=1,\ldots,d$. The characterization theorem of \cite{Sklar1959} implies that there exists a copula
function $\mathbb{ C}(\cdot)$ on $[0,1]^d$, such that
\begin{equation*}
\mathbb{ F}({\bf x}) = \mathbb{ C}( F_1(x_1),\ldots,F_d(x_d)),~~~\hbox{for}~~\mathbf{x}=(x_{1},\ldots,x_{d})\in \mathbb{R}^d.
\end{equation*}
This copula, which is unique, is a multivariate cumulative  distribution function whose univariate marginals are uniform on the interval $(0,1)$.
The Sklar's theorem provides the theoretical foundation for the
widespread use of the copula approach in generating multivariate distributions from univariate distributions.
The copula function pertaining to $\mathbb{F}(\cdot)$ may be defined by
\begin{equation}\label{copulrep}
\mathbb{ C}(\mathbf{ u}) = \mathbb{ F}(F^{-}_1(u_1),\ldots,
F^-_d (u_d))~\mbox{for }~\mathbf{ u}=(u_1,\ldots,u_d)\in[0, 1]^d,
\end{equation}
where, for $j = 1,\ldots,d$, $F^-_j (u) = \inf\{x\in \mathbb{R} :
F_j(x)\geq u\}$, with $u\in[0, 1]$, denotes the quantile function
of $F_j(\cdot)$. In the monographs by \cite{nelsen2006} and
\cite{Joe1997} the reader may find detailed ingredients of the
modelling theory as well as surveys of the commonly used copulas.
For in depth and overview historical notes we refer to
\cite{Schweizer1991}. Copulas have proved to be a  flexible and
versatile tool in the analysis of dependency structures. To be more
specific, copula $\mathbb{ C}(\cdot)$ ``couples'' the joint
distribution function $\mathbb{ F}(\cdot)$ to its univariate
marginals, capturing as such the dependence structure between the
components of $\mathbf{X }= (X_1,\ldots,X_d)$. Indeed, most
conventional measures of dependence can be explicitly expressed in
terms of the copula.
 This
feature has motivated successful applications in actuarial science
and survival analysis (see, e.g., \cite{frees1998}, \cite{Cui2004}).
In the literature on risk management and, more generally, in
mathematical economics and  mathematical finance modelling,
 a number of  illustrations are provided (refer to books of
\cite{Cherubini2004} and \cite{McNeil2005}), in particular,  in the context
of asset pricing and credit risk management.
First, we shall introduce some notations and definitions which will be used for the statement of our forthcoming results.
 Let  $\mathbf{X}_{i}= (X_{1i},\ldots,X_{di})$, $i= 1,2,\ldots$,  be independent random vectors with common distribution function
 $\mathbb{ F}(\cdot)$   whose   margins $F_{1}(\cdot),\ldots, F_d(\cdot)$ are   continuous  and whose   copula is denoted by
  $\mathbb{ C}(\cdot)$. Define $U_{ji} = F_j(X_{ji})$, for $i=1,\ldots,n$ and $j=1,\ldots,d$. The
   random  vectors $\mathbf{U}_i=(\xi_{1i},\ldots,\xi_{di})$   constitute  an  i.i.d.  sample  from $\mathbb{ C}(\cdot)$.
 Setting $\mathds{1}\{\cdots\}$ for the
indicator function of $\{\cdots\}$, we define, for each $n\geq1$, the following empirical distribution functions,
 for $\mathbf{ x} \in \mathbb{R}^{d}$ and for $\mathbf{ u}\in[0, 1]^d$,
\begin{eqnarray}
\mathbb{ F}_n(\mathbf{x})&:=&\frac{1}{n}\sum_{i=1}^n\prod_{j=1}^{d}\mathds{1}\{X_{ji}\leq
x_{j}\},\\
F_{jn}(x_{j})&:=&\frac{1}{n}\sum_{i=1}^n\mathds{1}\{X_{ji}\leq x_{j}\},\\
\mathbb{G}_n(\mathbf{u})&:=&\frac{1}{n}\sum_{i=1}^n\prod_{j=1}^{d}\mathds{1}\{U_{ji}\leq
u_{j}\},\\
G_{jn}(u_{j})&:=&\frac{1}{n}\sum_{i=1}^n\mathds{1}\{U_{ji}\leq u_{j}\}.
\end{eqnarray}
The marginal quantile functions associated to $F_{jn}(\cdot)$ and $G_{jn}(\cdot)$, for $j=1,\ldots,d$ and $u_{j}\in[0,1]$, are defined by
\begin{eqnarray}
\label{marges empirique1} F_{jn}^{-}(u_{j})&:=& \inf \{x\in\mathbb{R}:
F_{jn}(x)\geq u_j \},
\end{eqnarray}
\begin{eqnarray}
\label{marges empirique2} G_{jn}^{-}(u_j)&:=& \inf \{u\in[0,1]:
G_{jn}(u)\geq u_j \}.
\end{eqnarray}
According to  \cite{deheuvels1979a},
in view of the characterization (\ref{copulrep}), we define an empirical copula function of
$\mathbb{ F}_n(\cdot)$, based upon $\mathbf{ X}_{1}\ldots\mathbf{ X}_{n}$,  as any
copulas $\mathbb{ C}_{n}(\cdot)$ fulfilling the fundamental identity
\begin{equation}\label{copule empirique}
 \mathbb{ C}_{n}({\bf{u}}) := \mathbb{ F}_{n}(F_{1n}^-(u_1),\ldots,F_{jn}^-(u_d)),~~\mbox{for}~~\mathbf{u}\in [0,1]^d.
\end{equation}
It may be assumed without loss of generality that the marginal distributions of
 $\mathbb{F}(\cdot)$ are uniform on the interval $[0, 1]$, or equivalently that $\mathbb{F}(\cdot)$ is a copula, i.e.,
  we can work, in the sequel, directly with the sample
$\mathbf{ U}_{1},\ldots,\mathbf{ U}_n$ from $\mathbb{ C}(\cdot)$.
Here,  we may refer to
\cite{Sklar1959,Sklar1973}, \cite{Philipp1980}, \cite{wicura1973} and \cite{Moore1975} among others.
It follows that the empirical
copula in equation (\ref{copule empirique}
) is given by
\begin{equation}\label{copule empirique2}
 \mathbb{ C}_{n}({\bf{u}}) := \mathbb{G}_{n}(G_{1n}^-(u_1),\ldots,G_{jn}^-(u_d)),~~\mbox{for}~~\mathbf{u}\in [0,1]^d.
\end{equation}
The empirical copula function $\mathbb{ C}_{n}(\cdot)$ was briefly discussed by \cite{Ruymgaart1973}, pp. 6--13,
in the introduction of his
doctoral thesis. The asymptotic behavior of $\mathbb{ C}_{n}(\cdot)$ was studied in several papers, including
 \cite{deheuvels1979a}, \cite{stute1984}, \cite{stute1987}, \cite{Ruch1,Rush2} or \cite{tsukahara2005} and the
 references therein.
We may now define the \textit{normalized empirical copula process} $ \mathbb{A}_n(\cdot)$, according to
\cite{ruchendorf2009}, by setting
\begin{equation}\label{processuscopule}
  \mathbb{A}_n(\mathbf{u}):= n^{1/2}(\mathbb{ C}_n(\mathbf{u})-\mathbb{ C}(\mathbf{u})),~~
\mbox{ for }~~ \mathbf{u} \in [0,1]^d.
\end{equation}
The asymptotic behavior of the process $\{\mathbb{A}_n({\bf u}):{\bf u}\in [0,1]^d; n > 0\}$ has been
investigated extensively in stochastic literature. \cite{Deheuvels1981b,deheuvels1981} obtained the exact law and the
limiting process of $\{\mathbb{A}_n({\bf u}):{\bf u}\in [0,1]^d; n > 0\}$ under independence assumption
of margins. \cite{Ruch1,Rush2} and \cite{stute1987} proved weak convergence of the  process
$\{\mathbb{A}_n({\bf u}):{\bf u}\in [0,1]^d; n > 0\}$ in the space $D([0,1]^2),$
where   the space   of  c\`adl\`ag functions $D([0,1]^2)$ equipped with the Skorohod topology.
\cite{Wellner1996} established weak convergence in the space $\ell^{\infty}([a,b]^2),$ when $0< a< b<1,$
 under restrictions on the distribution functions.
\cite{fermanianradulovicdragan2004} showed  that the weak convergence of $\mathbb{A}_n(\cdot)$ to a centered Gaussian process
holds on $\ell^{\infty}([0,1]^2),$ when $\mathbb{ C}(\cdot)$ has continuous partial derivatives on $[0,1]^2.$
Recently, \cite{Segers2010} showed that the weak convergence of the \emph{normalized empirical copula process}
under the assumption that the first-order partial derivatives of the copula exist and are continuous on
certain subsets of the unit hyper-cube. We can say that the strong approximation holds for the process $\{\mathbb{A}_n({\bf u}):{\bf u}\in [0,1]^d; n > 0\}$ with rate $(b_{n})$, this means that, on a suitable probability space $(\Omega,\mathscr{A}, \mathbb{P})$,
\begin{equation}\label{invariance}
\sup_{\mathbf{u}\in [0,1]^{d}}|\mathbb{A}_{n}(\mathbf{u})-\mathbb{Z}_{n}(\mathbf{ u})|=O(b_{n}),~a.s.,
\end{equation}
 where $\mathbb{Z}_{n}(\cdot)$ is a sequence of Gaussian processes and  $b_{n}\rightarrow 0$ is a
  deterministic rate.  The strong approximations are quite useful and have received considerable
  attention in probability theory. Indeed, many well-known probability theorems can be considered as
consequences of results about strong approximation of sequences of sums by corresponding
Gaussian sequences. We shall mention that the rates of convergence for the distribution of \emph{smooth}
functionals of $\{\mathbb{A}_n({\bf u}):{\bf u}\in [0,1]^d; n > 0\}$ can  also be deduced from the
approximation in (\ref{invariance}).
The approximation by Kiefer processes is of particular interest, since any kind of law of the
iterated logarithm which holds for the partial sums of Gaussian processes may then be transferred
to the empirical processes $\{\mathbb{A}_n({\bf u}):{\bf u}\in [0,1]^d; n > 0\}$.
We refer to \cite{KMT}, \cite[Chapter 12]{Dasgupta2008}, \cite[Chapter 3]{Csorgoho1993},
\cite[Chapters 4-5]{Csorgo1981REVESZ} and  \cite[Chapter 12]{ShorackGalen1986s}
for expositions and references about this problem. We refer to \cite{Hall1983S} for a
survey of some applications of the strong approximation and many references. The interested
reader is referred to \cite{Dehe100000} and the references therein concerning the strong
approximations for the process  $\{\mathbb{A}_n({\bf u}):{\bf u}\in [0,1]^d; n > 0\}$.
In the last reference, a full characterization of empirical copula in general framework is provided.
 There is a
huge literature on the strong approximations and their applications. It is not the purpose
 of this paper to survey this extensive literature.\vskip7pt
\noindent  The remainder of the present paper is organized as follows.
In section \ref{sect22}, we will introduce notations and definitions
regarding some Gaussian processes, which play a central role in the
strong approximations theory.
In section \ref{sect22ZZZ}, we will give our main result
concerning the strong approximations of empirical copula processes by a single
Gaussian process, which is stated in Theorem \ref{th} below. In
section \ref{sect2}, we will give some applications of Theorem
\ref{th}. More precisely, we transfer our result to smoothed version of $\{\mathbb{A}_n({\bf u}):{\bf u}\in [0,1]^d; n > 0\}$ as well as to the
law of the iterated logarithm for the normalized empirical copula process.
To avoid interrupting the flow of the presentation, all mathematical developments are postponed to Section \ref{appe}.
\section{Main results}\label{sect}
\subsection{Gaussian Processes}\label{sect22}
Let $\mathbb{ C}(\cdot)$ be any copula. The $d$-variate Wiener process
$\{\mathbb{W}_\mathbb{ C}(\mathbf{ u}):\mathbf{ u}\in[0,1]^{d}\}$ on the
unit cube of $\mathbb{R}^{d}$ associated with the copula function $\mathbb{ C}(\cdot)$, i.e.,
$\mathbb{W}_{\mathbb{ C}}(\cdot)$ is a $d$-variate Gaussian process on $[0,1]^{d}$ with
\begin{eqnarray*}
\mathbb{E}(\mathbb{W}_\mathbb{ C}({\bf u}))=
0,\\
\mathbb{E}(\mathbb{W}_\mathbb{ C}({\bf
u})\mathbb{W}_\mathbb{ C}({\bf v}))&=& \mathbb{ C}({\bf u}\wedge{\bf
v}),
\end{eqnarray*}
where ${\bf u}\wedge{\bf v}:= (u_1\wedge v_1,\ldots,u_d\wedge
v_d)$ for ${\bf u}=(u_1,\ldots,u_d)\in [0,1]^d$ and ${\bf
v}=(v_1,\ldots,v_d)\in [0,1]^d,$
 and $\mathbb{W}_\mathbb{ C}(u_{1},\ldots,u_{d})=0$ whenever $u_{j}=0$, $j=1,\ldots,d$. \\
A $d$-variate Brownian bridge process on $[0,1]^d$ associated with the copula function $\mathbb{ C}(\cdot)$ is
defined, in terms of $\mathbb{W}_\mathbb{ C}(\cdot)$, by setting
\begin{equation}
\mathbf{B}_\mathbb{ C}({\bf u}):= \mathbb{W}_\mathbb{ C}({\bf
u})-\mathbb{ C}({\bf u})\mathbb{W}_\mathbb{ C}({\bf 1}),~~\hbox{for}~~{\bf u} \in [0,1]^d,
\end{equation}
where ${\bf 1}:=(1,\ldots,1).$ This
process has continuous sample paths and fulfills
\begin{eqnarray*}
\mathbb{E}(\mathbf{B}_\mathbb{ C}({\bf u}))&=&
0,\\\mathbb{E}(\mathbf{B}_\mathbb{ C}({\bf
u})\mathbf{B}_\mathbb{ C}({\bf v}))&=& \mathbb{ C}({\bf u}\wedge{\bf
v})-\mathbb{ C}({\bf u})\mathbb{ C}({\bf v}),~~~\mbox{for}~~{\bf
u},{\bf v} \in [0,1]^d.
\end{eqnarray*}The interested reader may refer to \cite{Piterbarg1996} and \cite{Adler1990}
for details on the Gaussian processes mentioned above. To state our result we need to define the
Kiefer random field. Consider a ($d+1$)-variate Gaussian process
 $\mathbb{W}_{\mathbb{ C}}(\mathbf{ u},z)$ on $[0,1]^{d}\times [0,\infty)$ such that
 $\mathbb{W}_{\mathbb{ C}}(\mathbf{ u},z)=0$ whenever any of $u_{1},\ldots,u_{d}$ or $z$ is zero. This process has continuous sample paths and fulfills
\begin{eqnarray*}
\mathbb{E}\left(\mathbb{W}_{\mathbb{ C}}(\mathbf{ u},z)\right)&=&0,\\
\mathbb{E}\left(\mathbb{W}_{\mathbb{ C}}(\mathbf{ u},z)\mathbb{W}_{\mathbb{ C}}(\mathbf{ v},t)\right)&=& \min(z,t)\mathbb{ C}(\mathbf{ u}\wedge \mathbf{ v}).
\end{eqnarray*}
 A ($d+1$)-variate Kiefer process $\mathscr{K}_{\mathbb{ C}}(\cdot,\cdot)$ on
$[0,1]^d\times[0,\infty)$ associated with the copula function
$\mathbb{ C}(\cdot)$, is defined, in term of $\mathbb{W}_{\mathbb{ C}}(\cdot,\cdot)$, by setting
\begin{equation}\label{ki212}
 \mathscr{K}_{\mathbb{ C}}({\bf{u}},t):= \mathbb{W}_\mathbb{ C}({\bf
u},t)-\mathbb{ C}({\bf u})\mathbb{W}_\mathbb{ C}({\bf 1},t)
\end{equation}
and fulfills
\begin{eqnarray*}
\mathbb{E}\left(\mathscr{K}_{\mathbb{ C}}({\bf{u}},z)\right)&=&0,\\\mathbb{E}\left(\mathscr{K}_{\mathbb{ C}}({\bf{u}},z)\mathscr{K}_{\mathbb{ C}}({\bf{v}},t)\right)&=&
(z\wedge t) \left\{\mathbb{ C}({\bf{u}}\wedge {\bf{v}}) -
\mathbb{ C}({\bf{u}})\mathbb{ C}({\bf{v}})\right\},
\end{eqnarray*}
 for
${\bf{u}},{\bf{v}} \in[0,1]^d$  and $z,~t\geq 0.$
We recall the distributional identity, for all fixed $z\geq 0,$
$$z^{-1/2}\mathscr{K}_{\mathbb{ C}}({\bf{u}},z)\stackrel{\mathscr{D}}{=}\mathbf{B}_\mathbb{ C}({\bf{u}}).$$  For each $n>0$, $u_j \in [0,1]$ and
$j=1,\ldots,d$, the \emph{copula Gaussian process} is defined by
 \begin{eqnarray} \label{k^*}
 \nonumber
 \mathscr{K}^{*}_{\mathbb{ C}}({\bf{u}},n) &:=&
\mathscr{K}_{\mathbb{ C}}({\bf{u}},n)-\sum_{j=1}^{d}\mathscr{K}_{\mathbb{ C}}(1,\ldots,1,u_{j},1,\ldots,1 ,n) \frac{\partial \mathbb{ C}({\bf
u})}{\partial{u_{j}}} \\
&=:&
\mathscr{K}_{\mathbb{ C}}({\bf{u}},n)-\sum_{j=1}^{d}\mathscr{K}_{\mathbb{ C}}^{(j)}({\bf
1 },u_{j},{\bf 1 },n) \frac{\partial \mathbb{ C}({\bf
u})}{\partial{u_{j}}}.
 \end{eqnarray}
We are now in position to state our main results.
\subsection{Strong approximation results}\label{sect22ZZZ}
In the sequel, the precise meaning of ``suitable probability
space'' is that an independent sequence of Wiener processes, which is independent
of the originally given sequence of i.i.d. r.v., can be constructed on
the assumed probability space. This is a technical requirement which allows
for the construction of the Gaussian processes in our theorems, and is not restrictive since
one can expand the probability space to make it rich enough (see e.g., Appendix 2 in \cite{Csorgoho1993}).
\vskip7pt
\noindent The main result to be proved here may now be stated precisely as
follows.
\begin{theorem} \label{th}
Assume that $\mathbb{ C}(\cdot)$, associated with $\mathbb{ F}(\cdot)$, is twice continuously differentiable on $(0, 1)^d$ and all the
partial derivatives of second order are continuous on $[0,1]^d$.
On a suitable probability space, we may define the empirical copula processes $\{\mathbb{A}_n({\bf u}):{\bf u}\in [0,1]^d, n \geq1\}$
 in combination with the Gaussian process
$\{\mathscr{K}^*_{\mathbb{ C}}({\bf u},t):{\bf u}\in [0,1]^d, t \geq 0 \}$, in such a way that, almost surely
as $n \rightarrow \infty$
\begin{equation}\label{limitdistri}
\sup_{{\bf {u}}\in[0,1]^d}\left|\sqrt{n}\mathbb{A}_n({\bf
u})-\mathscr{K}^*_{\mathbb{ C}}({\bf u},n)\right|= O \left(n^{1/2-
1/(4d)}(\log n)^{3/2}\right),
\end{equation}
where $\mathscr{K}^*_{\mathbb{ C}}({\bf u},t)$ is defined in (\ref{k^*}).
\end{theorem}
The proof of Theorem \ref{th} is captured in the forthcoming Section \ref{appe}.
\begin{remark}
In the particular case of independence, i.e.,  $$\mathbb{ C}({\bf u})=\prod_{j=1}^du_j,$$ the process
$\{\mathscr{K}^*_{\mathbb{ C}}({\bf u},n):{\bf u}\in [0,1]^d; n \geq 0 \}$
is equal to
 \begin{equation}
 \nonumber
 \mathscr{K}^{*}_{\mathbb{ C}}({\bf{u}},n)
=
\mathscr{K}_{\mathbb{ C}}({\bf{u}},n)-\sum_{j=1}^{d}\mathscr{K}_{\mathbb{ C}}^{(j)}({\bf
1 },u_{j},{\bf 1 },n) \prod_{i\neq j}^du_i, ~~~~{\bf u}\in[0,1]^d,
 \end{equation}
 with mean zero and covariance function
\begin{equation*}
\mathbb{E}\left(\mathscr{K}^{*}_{\mathbb{ C}}({\bf{u}},s)\mathscr{K}^{*}_{\mathbb{ C}}({\bf{v}},t)\right)=
(s\wedge t)\left\{\prod_{j=1}^d(u_j\wedge v_j)+(d-1)\prod_{j=1}^du_jv_j-\sum_{j=1}^d(u_j\wedge v_j)\prod_{i\neq j}u_iv_j\right\},
\end{equation*}
where   ${\bf{u}},{\bf{v}}\in [0,1]^d$ and $s,~t\geq 0$. For more details the reader may refer to \cite{csorgo1979}.
Note that in the case where $\{\mathbb{A}_n({\bf u}):{\bf u}\in [0,1]^d; n > 0\}$ is generated by a sample of random vectors
with independent marginals then the limit distribution in (\ref{limitdistri}) is free.
\end{remark}
\begin{remark}
 Theorem
\ref{th} may be used to derive the limiting laws of some
statistics like Kendall's sample rank correlation coefficient and
Spearman's sample rank correlation coefficient. More generally, let us define, for any function $J(\cdot)$ on $ [0,1]^3$
\begin{equation*}
\mathbb{S}(\mathbb{ C}):=\int_0^1\int_0^1 J(u,v,\mathbb{ C}(u,v))dudv.
\end{equation*}
The corresponding sample quantity $\mathbb{S}(\mathbb{ C}_n)$ may
be called Spearman type rank statistic, the interested reader may
refer to \cite{stute1987} and \cite{Tsukahara2000} for more
details. To be more precise, suppose that $z\rightarrow J(u,v,z)$
has a continuous derivative $J^3(u,v,z)$ with
$\sup_{u,v,z}|J^3(u,v,z)|=\sup_{u,v,z}|\partial J(u,v,z)/\partial
z|<\infty$. Then we can write
\begin{eqnarray*}
\lefteqn{\sqrt{n}(\mathbb{S}(\mathbb{ C}_n)-\mathbb{S}(\mathbb{ C}))}\\&=&\sqrt{n}\left(\int_0^1\int_0^1
J(u,v,\mathbb{ C}_n(u,v))dudv-\int_0^1\int_0^1 J(u,v,\mathbb{ C}(u,v))dudv\right)\\
&=&\int_0^1\int_0^1 J^3(u,v,\delta_n(u,v))\mathbb{A}_n(u,v)dudv,
\end{eqnarray*}
where $\delta_n(u,v)$ is a point between $\mathbb{ C}_n(u,v)$ and $\mathbb{ C}(u,v)$, so that $\delta_n$ converge to $ \mathbb{ C}$ uniformly with probability one. Making use of  Theorem
\ref{th} we have
\begin{equation*}
\left|\sqrt{n}(\mathbb{S}(\mathbb{ C}_n)-\mathbb{S}(\mathbb{ C}))-\int_0^1\int_0^1 J^3(u,v,\mathbb{ C}(u,v))
\frac{1}{\sqrt n} \mathscr{K}^{*}_{\mathbb{ C}}(u,v,n)dudv\right|=o_{\mathbb{P}}(1).
\end{equation*}
We put, for any function $J(\cdot)$ on $[0,1]^3$,
\begin{equation*}
\mathbb{T}(\mathbb{ C}):=\int_0^1\int_0^1 J(u,v,\mathbb{ C}(u,v))d\mathbb{ C}(u,v),
\end{equation*}
the integration being understood as multiple integral based on the bivariate  copula.
We call $\mathbb{T}(\mathbb{ C}_n)$ a Kendall type rank statistic. Similarly, using  Theorem
\ref{th} we can obtain the limiting law of $\sqrt{n}(\mathbb{T}(\mathbb{ C}_n)-\mathbb{T}(\mathbb{ C}))$.
\end{remark}
\begin{remark}
Theorem \ref{th} may be used to derive the limiting law of some smooth functionals.
We can see this, in particular, for the Kolmogorov-Smirnov and Cram\'er-von Mises type statistics, respectively,   given by
$$
\sup_{\mathbf{ u}\in [0,1]^{d}}|\mathbb{A}_n({\bf
u})|, ~~\mbox{and}~~\int_{[0,1]^{d}}\mathbb{A}_n^{2}({\bf
u})d\mathbf{ u}.
$$
We get
\begin{equation}
\left|\sup_{\mathbf{ u}\in [0,1]^{d}}|\mathbb{A}_n({\bf
u})|-\sup_{\mathbf{ u}\in [0,1]^{d}}\frac{1}{n^{1/2}}|\mathscr{K}^*_{\mathbb{ C}}({\bf u},n)|\right|
= O \left(n^{-
1/(4d)}(\log n)^{3/2}\right),
\end{equation}
and
\begin{equation}\label{blz}
\left|\int_{[0,1]^{d}}\mathbb{A}_n^{2}({\bf
u})d\mathbf{ u}-\frac{1}{n}\int_{[0,1]^{d}}\mathscr{K}^{*2}_{\mathbb{ C}}({\bf u},n)d\mathbf{ u}\right|= O \left(n^{-
1/(4d)}(\log n)^{3/2}(\log \log n)^{1/2}\right).
\end{equation}
At this point, we mention that the proof of (\ref{blz}) closely follows  the lines of 
\cite{Bouzebda2011MMS} and \cite{Bouzebda2010} among others. Therefore, we omit the details.
Another interesting application of the approximation
of $\mathbb{A}_n(\mathbf{u})$ in terms of Gaussian process in both $\mathbf{ u}$ and $n$ is the change-point problem, as in \cite{Remillard2010} and \cite{Cosorgo1997}, and leaves this study open for future research.
\end{remark}
Note that the covariance structure of  the process
$\{\mathscr{K}^{*}_{\mathbb{ C}}(\mathbf{u},n) : \mathbf{u} \in [0, 1]^d,~~n\geq 1\}$ depends on the first
 derivatives of the copula $\mathbb{ C}(\cdot)$ which are, in general, unknown in practice.
 To circumvent this problem, one can use  a multiplier central limit theorem, please refer to \cite{Wellner1996},
 as suggested in \cite{Remillard2009}, \cite{Scaillet2005} and recently \cite{Remillard2010}.
We mention that the usual bootstrap based on resampling was proposed in \cite{fermanianradulovicdragan2004}.
 Here, for easy reference and completeness, we recall the procedure given in \cite{Remillard2009}, which is
 more appropriate for our setting. Let $N$ be a large integer, and let $Z^{(k)}_{i}$, $i=1,\ldots,n$, $k= 1,\ldots,N$,
 be i.i.d. random variables with mean $0$ and variance $1$, independent of the data $\mathbf{ U}_{1},\ldots,\mathbf{ U}_{n}$.
 Moreover, for any $k= 1,\ldots,N$, let
\begin{eqnarray*}
\boldsymbol{\alpha}_n^{(k)}(\mathbf{u})&:=&\frac{1}{\sqrt{n}}\sum_{i=1}^{n}Z^{(k)}_{i}\left\{\mathds{1}
\left\{G_{1n}(U_{1i})\leq u_1,\ldots, G_{dn}(U_{di})\leq u_d\right\}-\overline{\mathbf{ \mathbb{ C}}}_{n}
({\bf u})\right\}\\&=&\frac{1}{\sqrt{n}}\sum_{k=1}^{n}\left(Z^{(k)}_{i}-\overline{\varepsilon}_n\right)
\mathds{1}\left\{G_{1n}(U_{1i})\leq u_1,\ldots, G_{dn}(U_{di})\leq u_d\right\},
\end{eqnarray*}
where
$$ \overline{\mathbf{ \mathbb{ C}}}_{n}({\bf u}): = \frac{1}{n} \sum_{i=1}^{n}\mathds{1}\left\{G_{1n}
(U_{1i})\leq u_1,\ldots, G_{dn}(U_{di})\leq u_d\right\},$$
and  $\overline{\varepsilon}_n=\frac{1}{n}\sum_{i=1}^n\varepsilon_i$, and for any $j=1,\ldots,d$,
\begin{eqnarray*}
\alpha_{jn}^{(k)}(u_j)&:=&\boldsymbol{\alpha}_n^{(k)}(1,\ldots,1,u_j,1,\ldots,1)\\&=&\frac{1}{\sqrt{n}}
\sum_{i=1}^{n}\left(Z^{(k)}_{i}-\overline{\varepsilon}_n\right)\mathds{1}\left\{G_{jn}(U_{ji})\leq u_j\right\}.
\end{eqnarray*}
Finally, for all $\mathbf{ u}\in[0,1]^d$, and for all $k=1,\ldots,N$, let
\begin{equation}
\mathbb{A}_{n}^{(k)}({\bf
u}):=\boldsymbol{\alpha}_n^{(k)}(\mathbf{u})-\sum_{j=1}^{d}\alpha_{jn}^{(k)}(u_j)\mathbb{C}^{[j]}_{n}(\mathbf{ u}),
\end{equation}
where $\mathbb{C}^{[j]}_{n}(\mathbf{ u})$ is a consistent estimator of the partial derivative $\frac{\partial \mathbb{ C}({\bf
u})}{\partial{u_{j}}}$, for $j=1,\ldots,d$, and given, as in \cite{Remillard2009}, by
$$
\mathbb{C}^{[j]}_{n}(\mathbf{ u}):=\frac{\mathbb{C}_{n}(\mathbf{ u}+he_{j})-\mathbb{C}_{n}(\mathbf{ u}-he_{j})}{2h},
$$
where $e_{j}$ is the $j$-th column of the $d \times d$ identity matrix. From Theorem 2.1 in
\cite{Remillard2009}, a practical choice of  $h$  is $n^{-1/2}$. We will not investigate the
 question of the strong approximations of the processes $\mathbb{A}_{n}^{(k)}(\cdot)$  in the present paper.
\begin{remark}
In \cite{Dette20102},
a new procedure, which has the most attractive theoretical properties, has been proposed
to circumvent the problem of the estimation of the derivatives of the unknown copula.
In the same paper, the finite-sample properties of some methods are compared in a simulation study,
and the multiplier approach, by  \cite{Remillard2009}, yields the best results in most cases.
\end{remark}
Recently, \cite{Remillard2009} considered the two sample problem in the copula setting.
Let
$\mathbb{C}_n(\cdot)$ and $\mathbb{D}_m(\cdot)$ denote the
empirical copulas functions based on independent samples of sizes
$n$ and $m$, respectively. Theorem \ref{th} may be used for statistical
comparison procedures of the unknown copulas $\mathbb{C}(\cdot)$
and $\mathbb{D}(\cdot)$ based on $\mathbb{C}_n(\cdot)$ and
$\mathbb{D}_m(\cdot)$.
Consider the empirical process  $\mathbb{A}_{n;m}(\cdot)$, defined by
\begin{equation}\label{emxi}
\mathbb{A}_{n;m}({\bf u}):= \left\{\frac{nm}{n+m}\right\}^{1/2}
\left\{\mathbb{C}_{n}({\bf u})-\mathbb{D}_{m}({\bf u})
\right\},~~\mbox{for}~~{\bf u}\in [0,1]^d.
\end{equation}
In order to test the null hypothesis $\mathscr{H}_0:\mathbb{C}=\mathbb{D}$, we use Cram\'er-von Mises type statistic, given by
\begin{equation}\label{CviC}
\Omega_{n;m} :=\int_{[0,1]^d}\{\mathbb{A}_{n;m}({\bf u})\}^2 ~d{\bf
u}.
\end{equation}
We need to define the following Gaussian process
\begin{eqnarray}
\mathscr{K}_{n;m}^{*}({\bf u}):= [m/(n+m)]^{1/2}\frac{1}{n^{1/2}}\mathscr{K}^*_{\mathbb{ C}}({\bf u},n)+[n/(n+m)]^{1/2}
\frac{1}{m^{1/2}}\mathscr{K}^*_{\mathbb{D}}({\bf u},m).
\end{eqnarray}
Using Theorem \ref{th}, one can show,  as $\min(n,m)\rightarrow\infty$ and $n/(n + m)\rightarrow \lambda \in [0, 1]$, then we have
almost surely,
\begin{eqnarray*}
\lefteqn{\left|\int_{[0,1]^d}\{\mathbb{A}_{n;m}({\bf u})\}^2 ~d{\bf
u}-\int_{[0,1]^d}\{\mathscr{K}_{n;m}^{*}({\bf u})\}^2 ~d{\bf
u}\right|}\\&=&O\left(\max\left(\frac{(\log n)^{3/2}}{n^{1/(4d)}},\frac{(\log m)^{3/2}}{m^{1/(4d)}}\right)
\times\max\left((\log \log n)^{1/2},(\log \log m)^{1/2}\right)\right).
\end{eqnarray*}

\section{Applications}\label{sect2}
\subsection{Smoothed empirical copula processes}
\noindent A seemingly natural kernel-type estimator $\widehat{\mathbb{ C}}_n(\cdot)$ of $\mathbb{ C}(\cdot)$ would be
\begin{equation}\label{kern231}
\widehat{\mathbb{ C}}_n(\mathbf{u}):=\frac{1}{h}\int_{[0,1]^d}
k\left(\frac{\mathbf{u}-\mathbf{v}}{h^{1/d}}\right)\mathbb{ C}_n(\mathbf{v})d\mathbf{v},~~\mbox{
for }~~\mathbf{u}\in[0,1]^d,
\end{equation}
where $k(\cdot)$ is a kernel function and $h=h(n)$ is the smoothing parameter.
For notational convenience, we have chosen the same bandwidth sequence for each margins. This assumption can be dropped easily.
The kernel estimation of copula function is a
rich topic of researches, we only mention  \cite{Omelka2009}, \cite{Chen200777} and  \cite{fermanianradulovicdragan2004}, see their lists of
references for related studies. As in the previous section, we define the smoothed empirical
copulas process, for $n\geq 1$,  by
\begin{equation}
\widehat{\mathbb{A}}_n(\mathbf{u}):=\sqrt{n}\left(\widehat{\mathbb{ C}}_n(\mathbf{u})-\mathbb{ C}(\mathbf{u})\right),~\mbox{
for }~\mathbf{u}\in[0,1]^d.
\end{equation}
We will describe the asymptotic properties of the  smoothed empirical
copulas process $\{\widehat{\mathbb{A}}_n({\bf u}):{\bf u}\in [0,1]^d; n > 0\}$ under the following assumptions.
\begin{enumerate}
\item[(F.1)]  There exists a constant $0<\mathfrak{C}<\infty$ such that
\begin{equation*}
\sup_{\mathbf{u} \in[0,1]^d}\left|\frac{\partial ^s \mathbb{ C}(\mathbf{u})}{\partial^{j_1}u_1\ldots\partial^{j_d}u_d}
\right|\leq \mathfrak{C},~~j_1+\cdots+j_d=s.
\end{equation*}
\end{enumerate}
Suppose  that $\{h(n)\}_{n\geq1}$ is a sequence of positive constants
which satisfies the following condition.
\begin{enumerate}
\item[(C.1)] $h=h(n)\rightarrow 0,$ $nh\rightarrow \infty$ and
$\sqrt{n}h^{s/d}\rightarrow 0$ as  $n\rightarrow\infty$.
\end{enumerate}
The following conditions on the kernel function $k(\cdot)$ are assumed in our analysis.
\begin{enumerate}
\item[(C.2)] $k(\cdot)$ is a  continuous density function and compactly supported; \item[(C.3)] $k(\cdot)$ is of order $s$, i.e.,
\begin{eqnarray*}
&&\int_{\mathbb{R}^d}k(\mathbf{u})d\mathbf{u}=1,\\
&&\int_{\mathbb{R}^d}u_1^{j_1}\ldots u_d^{j_d}k(\mathbf{u})d\mathbf{u}=0,~~j_1,\ldots,j_d\geq 0,~~j_1+\cdots +j_d=1,\ldots,s-1,\\
&&\int_{\mathbb{R}^d}|u_1^{j_1}\ldots u_d^{j_d}|k(\mathbf{u})d\mathbf{u}<\infty,~~j_1,\ldots,j_d\geq 0,~~j_1+\cdots +j_d=s.
\end{eqnarray*}
\end{enumerate}
It is now possible to state the main theoretical result of this
section which provides the limiting behavior of the smoothed
empirical copulas process $\{\widehat{\mathbb{A}}_n({\bf u}):{\bf
u}\in [0,1]^d; n > 0\}$.
\begin{corollary} \label{th1} Assume that ${\rm (F.1)}$ and ${\rm (C.1)}$-${\rm(C.3)}$ hold.
Then, on a suitable probability space, we may define the smoothed empirical copula processes
$\{\widehat{\mathbb{A}}_n({\bf u}):{\bf u}\in [0,1]^d; n > 0\}$
in combination with  the Gaussian process
$\{\mathscr{K}^*_{\mathbb{ C}}({\bf u},t):{\bf u}\in [0,1]^d; t\geq  0 \}$, in such a way that,
as $n \rightarrow \infty$
\begin{equation}\label{smoooo}
\sup_{{\bf {u}}\in[0,1]^d}\left|\widehat{\mathbb{A}}_n({\bf
u})-\frac{1}{\sqrt{n}}\mathscr{K}^*_{\mathbb{ C}}({\bf u},n)\right|=o_{\mathbb{P}}(1).
\end{equation}
\end{corollary}
The proof of Corollary \ref{th1} is postponed until  Section \ref{appe}.\\
The result of Corollary \ref{th1} is motivated by the following remark.
\begin{remark}
 The empirical copula  provides a universal
way  for estimation purposes. Unfortunately, its discontinuous feature
induces some difficulties: the graphical representations of the
copula may not be satisfactory  from a visual and intuitive point of view.
Moreover, there is no unique choice for building the inverse
function of marginal functions.  Finally, since
the empirical copula estimator is not differentiable, it cannot, for example, be used  for optimization purposes.  Studies have shown that a smoothed estimator  may be preferable to
the sample estimator. First, smoothing reduces the random variation in the data,
resulting in a more efficient estimator. Second, smoothing gives a smooth curve that displays  some interesting features.
For more details on the subject we may refer to \cite{Chen200777}.
\end{remark}
\begin{remark}
\begin{enumerate}\item Corollary \ref{th1} remains valid when replacing the
condition that the kernel function $k(\cdot)$ having compact support in ${\rm(C.2)}$ by
another condition ${\rm(C.4)}$ which content is as follows
\begin{enumerate}
\item[{\rm(C.4)}] There exists a sequence of positive real numbers $a_n$ such
that $a_nh$ tends to zero when $n$ tends to infinity, and
$$\sqrt{n}\int_{\{\parallel \mathbf{v}  \parallel > a_n\}}|k(\mathbf{v})|d\mathbf{v}\rightarrow 0.$$
\end{enumerate}
\item Note that the conditions of Corollary \ref{th1} are grouped to
control the deviations between the \textit{normalized empirical copula process} $ \{\mathbb{A}_n({\bf u}):{\bf u}\in [0,1]^d; n > 0\}$ and the
smoothed empirical copula process $\{\widehat{\mathbb{A}}_n({\bf u}):{\bf u}\in [0,1]^d; n > 0\}$.
\end{enumerate}
\end{remark}
\subsection{The law of iterated logarithm for the normalized empirical copula process}
From Theorem \ref{th}, we have almost surely
\begin{equation}\label{lil}
\limsup_{n\rightarrow\infty}\left\{\left(\frac{n}{2\log \log
n}\right)^{1/2}\sup_{{\bf
{u}}\in[0,1]^d}|\mathbb{ C}_n(\mathbf{u})-\mathbb{ C}(\mathbf{u})|\right\}=\limsup_{n\rightarrow\infty}\frac{\sup_{{\bf
{u}}\in[0,1]^d}|\mathscr{K}^*_{\mathbb{ C}}({\bf u},n)|}{(2n\log \log
n)^{1/2}}.
\end{equation}
Note that (\ref{lil}) readily implies the following corollary, which is a straightforward consequence of Theorem \ref{th}.
\begin{corollary} \label{th2} Under the same conditions of Theorem \ref{th}, we
have
\begin{equation}
\limsup_{n\rightarrow \infty}\left\{\left(\frac{n}{2\log \log
n}\right)^{1/2}\sup_{{\bf
{u}}\in[0,1]^d}|\mathbb{ C}_n(\mathbf{u})-\mathbb{ C}(\mathbf{u})|\right\}=\rho,\hbox{  } a.s.,
\end{equation}
where
\begin{equation*}
\rho^2:=\sup_{\mathbf{u}\in[0,1]^d}{\rm Var}\left(\mathscr{K}^*_{\mathbb{ C}}({\bf u},1)\right).
\end{equation*}
\end{corollary}
\begin{remark}
A result similar to Corollary \ref{th2} was obtained by \cite{deheuvels1979a} (refer to Theorem 3.1) using a different method.
\end{remark}
\begin{remark}
Statistics of the form
\[\mathbf{R}_n:=\frac{1}{n}\sum_{i=1}^nJ\left(G_{1n}(X_{1i}),\ldots,G_{dn}(X_{di})\right),\] belong to the general
class of \emph{multivariate rank statistics}. Their asymptotic
properties have been investigated at length by a number of
authors, among whom we may quote \cite{Ruymgaart_Shorack_Zwet1972},
\cite{Ruch1,Rush2}. In particular,
the previous authors have provided regularity conditions, imposed
on $J(\cdot)$, which imply the asymptotic normality of
$\mathbf{R}_n$. It is easy to see that
$$\mathbf{R}_n=\int_{[0,1]^{d}}J(\mathbf{ u})d\overline{\mathbb{ C}}_{n}(\mathbf{ u}).$$
Since the difference between $\overline{\mathbf{ \mathbb{ C}}}_{n}(\cdot)$ and $\mathbb{ C}_{n}(\cdot)$ is negligible,
 see \cite{fermanianradulovicdragan2004} or \cite{Dehe100000},
the asymptotic normality of~$\mathbf{R}_n$ can be established under the
weakest set of assumptions (see, Theorem 6 in \cite{fermanianradulovicdragan2004}) using  Theorem \ref{th}.
\end{remark}

\section{Proofs}\label{appe}
This section is devoted to the proofs of our results.
\subsection*{Proof of Theorem \ref{th}.}
Consider the empirical processes defined,
respectively, for $n \geq 1$, $\mathbf{u}\in[0,1]^d$ and $0 \leq
u_j\leq 1$, for $j=1,\ldots,d$, by
\begin{eqnarray}\label{empiri}
 \boldsymbol{\alpha}_n(\mathbf{u})&:=&
n^{1/2}(\mathbb{G}_n(\mathbf{u})-\mathbb{ C}(\mathbf{u})),\\
 \alpha_{jn}(u_{j})&:=& n^{1/2}
\{G_{jn}(u_{j}) -u_{j}\},\\
\beta_{jn}(u_{j}) &:=& n^{1/2} \{G_{jn}^{-}(u_{j}) -u_{j}\}.
\end{eqnarray}
Keep in mind the
definition (\ref{processuscopule}) of $\mathbb{A}_n(\cdot)$.
The \textit{normalized empirical copula process} can be written, for $ \mathbf{u} \in
[0,1]^d$,  as follows
\begin{eqnarray} \label{3 parties de G_n}
\nonumber
  \lefteqn{\mathbb{A}_n(\mathbf{u})=
n^{1/2}(\mathbb{G}_n(G_{1,n}^-(u_1),\ldots,G_{d,n}^-(u_d))-\mathbb{ C}(u_1,\ldots,u_d))}\\
\nonumber &=& \boldsymbol{\alpha}_n\left(G_{1,n}^-(u_1),\ldots,G_{d,n}^-(u_d)
\right)+
n^{1/2}\left\{\mathbb{ C}\left(G_{1,n}^-(u_1),\ldots,G_{d,n}^-(u_d)\right)
\right.\\\nonumber&&~~~~~~~~~~~~~~- \left.\mathbb{ C}(u_1,\ldots,u_d)\right\}\\
\nonumber &=&
\boldsymbol{\alpha}_n\left(u_1+n^{-1/2}\beta_{1n}(u_{1}),\ldots,u_d+n^{-1/2}\beta_{dn}(u_{d})
\right) \\
\nonumber &&+
n^{1/2}\left\{\mathbb{ C}\left(u_1+n^{-1/2}\beta_{1n}(u_{1}),\ldots,u_d+n^{-1/2}\beta_{dn}(u_{d})\right)
-\mathbb{ C}(u_1,\ldots,u_d))\right\}\\
\nonumber &=& \boldsymbol{\alpha}_n({\bf u})+ \left\{ \boldsymbol{\alpha}_n({\bf u}+
n^{-1/2}{\bf \beta}_{n}({\bf u}))- \boldsymbol{\alpha}_n({\bf u})\right\}
\\\nonumber &&+ n^{1/2}\left\{ \mathbb{ C}({\bf u}+ n^{-1/2}{\bf
\beta}_{n}({\bf
u}))- \mathbb{ C}({\bf u})\right\}\\
&=& \boldsymbol{\alpha}_n({\bf u})+ \boldsymbol{\Delta}_1({\bf u}, n)+ \boldsymbol{\Delta}_2({\bf u}, n),
\end{eqnarray}
where $({\bf u}+ n^{-1/2}{\bf \beta}_{n}({\bf u}))=
(u_1+n^{-1/2}\beta_{1n}(u_{1}),\ldots,u_d+n^{-1/2}\beta_{dn}(u_{d}))$.
The decomposition (\ref{3 parties de G_n}) is the main key to our proof. We first compute the
 right side term $\boldsymbol{\Delta}_2(\cdot, n)$  of (\ref{3 parties de G_n}).
Under differentiability
assumption on $\mathbb{ C}(\cdot)$ and  by successive Taylor expansions, we readily obtain the
equality
\begin{eqnarray*}
\boldsymbol{\Delta}_2({\bf u}, n)&=&\sum_{j=1}^{d}\frac{\partial
\mathbb{ C}(\mathbf{u})}{\partial
u_j}\sqrt{n}(G_{nj}^-(u_j)-u_j)\\&&+\frac{\sqrt{n}}{2}\sum_{j=1}^{d}\sum_{k=1}^{d}\frac{\partial^2\mathbb{ C}(\mathbf{u^*})}{\partial
u_j\partial u_k}(G_{nj}^-(u_j)-u_j)(G_{nk}^-(u_k)-u_k),
\end{eqnarray*}
which holds for some point $\mathbf{u}^*$  in  the interior of the line segment
joining $(G_{n1}^-(u_1),\ldots,G_{nd}^-(u_d))$ and
$(u_1,\ldots,u_d)$. It follows from the definition of
$\boldsymbol{\alpha}_n(\cdot)$ in (\ref{empiri}), for $ u_j \in
[0,1],~j=1,\ldots,d,$ that
\begin{eqnarray*}
\sqrt{n}(G_{nj}^-(u_j)-u_j)&=& -\sqrt{n}\left(G_{nj}(G_{nj}^-(u_j))- G_{nj}^-(u_j)\right)+ \sqrt{n}\left(G_{nj}(G_{nj}^-(u_j))- u_j\right) \\
&=&
-\boldsymbol{\alpha}_n(\mathbf{1},G_{nj}^-(u_j),\mathbf{1})+\sqrt{n}(G_{nj}(
G_{nj}^-(u_j))-u_j).
\end{eqnarray*}
Using the fact, for $ u_j \in
[0,1],~j=1,\ldots,d,$ that $$\left|G_{nj}(G_{nj}^-(u_j))-u_j\right|\leq \frac{1}{n}$$ and
the \cite{sChung1949}'s law of the iterated logarithm, one
finds, almost surely,
\begin{eqnarray*}
\boldsymbol{\Delta}_2({\bf u}, n)&=&-\sum_{j=1}^{d}\frac{\partial
\mathbb{ C}(\mathbf{u})}{\partial
u_j}\boldsymbol{\alpha}_n(\mathbf{1},G_{nj}^-(u_j),\mathbf{1})+O(n^{-1/2}\log
\log n),
\end{eqnarray*}
uniformly in ${\bf u}\in[0,1]^{d}$. It is well known
from Stute's work [\cite{Stute1982},
p. 99], that we have, almost surely, for $n$ sufficiently large and $j=1,\ldots,d$,
\begin{equation*}
\sup_{u_j\in[0,1]}|\boldsymbol{\alpha}_n(\mathbf{1},G_{nj}^-(u_j),\mathbf{1})-\boldsymbol{\alpha}_n(\mathbf{1},u_j,\mathbf{1})|=O(n^{-1/4}(\log
n)^{1/2}(\log \log n)^{1/4}).
\end{equation*}
Then, it follows that uniformly in ${\bf u}\in[0,1]^{d}$, almost surely, for $n$
sufficiently large
\begin{eqnarray}\label{delta2222}
\boldsymbol{\Delta}_2({\bf u}, n)&=&-\sum_{j=1}^{d}\frac{\partial
\mathbb{ C}(\mathbf{u})}{\partial
u_j}\boldsymbol{\alpha}_n(\mathbf{1},u_j,\mathbf{1})+O(n^{-1/4}(\log
n)^{1/2}(\log \log n)^{1/4}),
\end{eqnarray}
as was observed by \cite{stute1984}, p. 371.
We next evaluate  the term $\boldsymbol{\Delta}_1(\cdot, n)$ in the right hand side of (\ref{3 parties de G_n}).
Recall that $\boldsymbol{\Delta}_1({\bf u}, n)$ is the difference between $\boldsymbol{\alpha}_n({\bf
u}+ n^{-1/2}{\bf \beta}_{n}({\bf u}))$ and $\boldsymbol{\alpha}_n({\bf u})$.
Let $w_{n}(\cdot)$ be the modulus of continuity of
$\boldsymbol{\alpha}_n(\cdot)$, that is
$$ w_{n}({\bf{a}}) :=
\sup\left\{ \alpha_{n}(A): A=\prod_{j=1}^d [u_j,v_j] \in
[0,1]^d,~~ \hbox{with}~~ \left|[u_j,v_j]\right|= v_j-u_j \leq a_j, \forall j=
1,\ldots,d \right\},
$$
where $\mathbf{ a}:=(a_{1},\ldots,a_{d})$.
We will make use of the following fact which is a particular case of Theorem 2.1, p. 367 of \cite{stute1984}.\\ \noindent {\bf Fact
1.} Let  $\{a_n\}_{n\geq 1}$
be a sequence in $(0,1)$ such as $a_n \downarrow 0$, as
$n\rightarrow \infty$, and
$$ i)  na_n^d\uparrow \infty,~~~~ii) na_n^d/ \log n \rightarrow \infty
,~~~~ iii) \log (1/a_n)/ \log \log n   \rightarrow \infty. $$ Then, we have, almost surely, $$ \lim_{n\rightarrow \infty}
 \left\{2 a_n^d \log (1/ a_n^d)\right\}^{-1/2} w_{n}(a_n,\dots,a_n)=
1.$$
 Once more, an application of the \cite{sChung1949} law of the iterated
logarithm shows that, for each $j=1,\ldots,d,$ almost surely,
\begin{equation}\label{chung}
\lim\sup_{n\rightarrow\infty} \left\{(\log\log n)^{-1/2}
\sup_{0\leq u_{j} \leq 1} |\beta_{jn}(u_j) |\right\} = 2^{-1/2}.
\end{equation}
In view of (\ref{chung}), we have almost surely, for all $j=1,\ldots,d$ and $n$ large enough,
\begin{eqnarray*}
\sup_{0\leq u_j \leq 1} \mid n^{-1/2}\beta_{jn}(u_j)\mid &\leq&
\frac{(\log \log n)^{1/2}}{n^{1/2}}\\&\leq& \frac{(\log
n)^{2/d}}{n^{1/d}}:= a_n,
\end{eqnarray*}
an application of  Fact 1  shows that,  as $n\rightarrow
\infty$, we have, almost surely,
\begin{equation}\label{delta111}
\sup_{\mathbf{u} \in [0,1]^d} \mid \boldsymbol{\Delta}_1({\bf u}, n)\mid \leq w_{n}(\mathbf{ a}_n)=
O\left( n^{-1/2} (\log n)^{3/2}\right),
\end{equation}
where $\mathbf{a}_{n}:=(a_{n},\ldots,a_{n})$.
The next fact,
due to \cite{Csorgo1988}, p. 102, provides a strong approximation result
appropriate to our need. Recall the definitions (\ref{ki212}) and (\ref{empiri}). \\
\noindent {\bf Fact 2.} On a suitable probability space $(\Omega,
\mathscr{A}, \mathbb{P})$,
 it is possible to define $\{\boldsymbol{\alpha}_n({\bf u}): {\bf u} \in [0,1]^d \}$, jointly  with
 the sequence of Gaussian processes $\{\mathscr{K}_{\mathbb{ C}}({\bf{u}},t): {\bf{u}} \in [0,1]^d, t\geq
0\}$, in such a way that,  as $n\rightarrow \infty$, almost surely,
\begin{equation}\label{l}
\sup_{{{\bf{u}} \in [0,1]^d}}\left|\sqrt{n}\boldsymbol{\alpha}_n({\bf{u}}) -
\mathscr{K}_{\mathbb{ C}}({\bf{u}},n)\right| = O \left(n^{1/2-
1/(4d)}(\log n)^{3/2}\right).
\end{equation}
In view of the above Fact,
by combining (\ref{delta2222}) and (\ref{delta111}) with the triangle inequality, we readily obtain
\begin{eqnarray*}
\lefteqn{\sup_{{\bf {u}}\in[0,1]^d}\left|\sqrt{n}
\mathbb{A}_n({\bf{u}})- \mathscr{K}^*_{\mathbb{ C}}({\bf u},n)\right| \leq
\sup_{{{\bf{u}} \in [0,1]^d}}\left|\sqrt{n}\boldsymbol{\alpha}_n({\bf{u}}) -
\mathscr{K}_{\mathbb{ C}}({\bf{u}},n)\right| + \sqrt{n} w_{n}(\mathbf{ a}_n)}\\
 &&+ \sum_{j=1}^d \left|\frac{\partial \mathbb{ C}({\bf u})}{\partial{u_{j}}}\right| \sup_{0\leq u_j \leq 1}
 \mid \sqrt{n}\boldsymbol{\alpha}_n(\mathbf{1},u_j,\mathbf{1})-\mathscr{K}_{\mathbb{ C}}^{(j)}({\bf 1 },u_{j},{\bf 1
 },n)\mid \\
 &&+
O(n^{1/4}(\log
n)^{1/2}(\log \log n)^{1/4}),\\
&\leq &O \left(n^{1/2- 1/(4d)}(\log n)^{3/2}\right) +
O\left((\log n)^{3/2}\right)\\
&&+ O\left( n^{1/2- 1/(4d)}(\log n)^{3/2}\right)
+O\left( n^{1/4} (\log n)^{1/2}(\log \log
n)^{1/4}\right)\\&=&O \left(n^{1/2-
1/(4d)}(\log n)^{3/2}\right).
\end{eqnarray*}
Note that we have used the fact that the first-order partial derivatives
of a copula are bounded (see Theorem 2.2.7 of \cite{nelsen2006}).
Then, we have almost surely, for all $n$ sufficiently large,
\begin{equation*}
\sup_{{\bf {u}}\in[0,1]^d}|\sqrt{n} \mathbb{A}_n({\bf{u}})-
\mathscr{K}^*_{\mathbb{ C}}({\bf u},n)|=O \left(n^{1/2-
1/(4d)}(\log n)^{3/2}\right)
\end{equation*}
and thus the proof of Theorem \ref{th} is completed.  \hfill$\Box$\\
\subsection*{Proof of Corollary \ref{th1}.}
We shall first study the behavior of the difference between
the \emph{normalized empirical copula process} $ \mathbb{A}_n(\cdot)$ and the
smoothed empirical copula process $ \widehat{\mathbb{A}}_n(\cdot)$. Recall the definition (\ref{kern231}). Notice that, for
each $\mathbf{ u}\in[0,1]^{d}$,
\begin{eqnarray}
\widehat{\mathbb{A}}_n(\mathbf{u})&=&\sqrt{n}\left(\widehat{\mathbb{ C}}_n(\mathbf{u})-\mathbb{ C}(\mathbf{u})\right)\nonumber\\
&=&\nonumber\sqrt{n}\left(\frac{1}{h}\int_{[0,1]^d}
k\left(\frac{\mathbf{u}-\mathbf{v}}{h^{1/d}}\right)
\mathbb{ C}_n(\mathbf{v})d\mathbf{v}-\mathbb{ C}(\mathbf{u})\right)\\\nonumber&=&\left(\frac{1}{h}\int_{[0,1]^d}
k\left(\frac{\mathbf{u}-\mathbf{v}}{h^{1/d}}\right)\sqrt{n}(\mathbb{ C}_n(\mathbf{v})-\mathbb{ C}(\mathbf{v}))\right)d\mathbf{v}
\\\nonumber&&+\sqrt{n}\left(\frac{1}{h}\int_{[0,1]^d}
k\left(\frac{\mathbf{u}-\mathbf{v}}{h^{1/d}}\right)\mathbb{ C}(\mathbf{v})d\mathbf{v}-\mathbb{ C}(\mathbf{u})\right)
\\\nonumber&=&\left(\frac{1}{h}\int_{[0,1]^d}
k\left(\frac{\mathbf{u}-\mathbf{v}}{h^{1/d}}\right)\mathbb{A}_n(\mathbf{v})\right)d\mathbf{v}\\&&+\sqrt{n}\left(\frac{1}{h}\int_{[0,1]^d}
k\left(\frac{\mathbf{u}-\mathbf{v}}{h^{1/d}}\right)\mathbb{ C}(\mathbf{v})d\mathbf{v}-\mathbb{ C}(\mathbf{u})\right).
\end{eqnarray}
We will make use of the following straightforward inequality
\begin{eqnarray}\label{decom}
\lefteqn{\sup_{{\mathbf{u}}\in[0,1]^d}|\widehat{\mathbb{A}}_n(\mathbf{u})-\mathbb{A}_n(\mathbf{u})|}\nonumber\\
&\leq&\sup_{{\mathbf{u}}\in[0,1]^d}\left|\int_{\prod_{i=1}^d\left[\frac{u_i-1}{h^{1/d}},\frac{u_i}{h^{1/d}}\right]}(\mathbb{A}_n(\mathbf{u}-h^{1/d}\mathbf{v})-\mathbb{A}_n(\mathbf{u}))
k(\mathbf{v})d\mathbf{v}\right|\nonumber\\&&+\sup_{{\mathbf{u}}\in[0,1]^d}|\mathbb{A}_n(\mathbf{u})|\left|\int_{\prod_{i=1}^d\left[\frac{u_i-1}{h^{1/d}},\frac{u_i}{h^{1/d}}\right]}
k(\mathbf{v})d\mathbf{v}-1\right|\nonumber\\&&+
\sqrt{n}\sup_{{\mathbf{u}}\in[0,1]^d}\left|\int_{\prod_{i=1}^d\left[\frac{u_i-1}{h^{1/d}},\frac{u_i}{h^{1/d}}\right]}(\mathbb{ C}(\mathbf{u}-h^{1/d}\mathbf{v})-\mathbb{ C}(\mathbf{u}))
k(\mathbf{v})d\mathbf{v}\right|\nonumber\\&&+\sqrt{n}\sup_{{\mathbf{u}}\in[0,1]^d}|\mathbb{ C}(\mathbf{u})|\left|\int_{\prod_{i=1}^d\left[\frac{u_i-1}{h^{1/d}},\frac{u_i}{h^{1/d}}\right]}
k(\mathbf{v})d\mathbf{v}-1\right|\nonumber\\&:=&\boldsymbol{\nabla}_{1;n}+\boldsymbol{\nabla}_{2;n}+\boldsymbol{\nabla}_{3;n}+\boldsymbol{\nabla}_{4;n}.
\end{eqnarray}
We first evaluate $\boldsymbol{\nabla}_{3;n}$ in the right side of (\ref{decom}).
Under conditions (F.1), (C.1)-(C.3) and applying a Taylor series expansion of
order $s$, we can see that
\begin{equation*}
\boldsymbol{\nabla}_{3;n}=\frac{h^{s/d}}{s!}
\sqrt{n}\sup_{{\mathbf{u}}\in[0,1]^d}\left|\int\sum_{j_1+\cdots+j_d=s}
 u_1^{j_1}\ldots u_d^{j_d}\frac{\partial^s\mathbb{ C}(\mathbf{u}-h_n\theta\mathbf{v})}{\partial
 u_1^{j_1}\ldots \partial u_d^{j_d}}k(\mathbf{v})d\mathbf{v}\right|,
\end{equation*}
where $\theta=(\theta_1,\ldots,\theta_d)$ and $0<\theta_j<1$, for $j=1,\ldots,d$. Thus, a straightforward application of Lebesgue
dominated convergence theorem gives
\begin{equation}
n^{-1/2}h^{-(s/d)}\boldsymbol{\nabla}_{3;n}=\frac{1}{k!}\sup_{{\mathbf{u}}\in[0,1]^d}
\left|\sum_{j_1+\cdots+j_d=s}\frac{\partial^s\mathbb{ C}(\mathbf{u})}{\partial
 u_1^{j_1}\ldots \partial u_d^{j_d}}\int
 u_1^{j_1}\ldots u_d^{j_d} k(\mathbf{v})d\mathbf{v}\right|.
\end{equation}
Then by condition (C.1) and (C.3), we conclude that, for all $n$ sufficiently large,
\begin{equation}\label{bia}
\boldsymbol{\nabla}_{3;n}=O(n^{1/2}h^{s/d})=o(1).
\end{equation}
 Making use of Theorem \ref{th} in connection with the almost sure continuity
of the Gaussian process $\{\mathscr{K}^*_{\mathbb{ C}}({\bf u},t):{\bf u}\in [0,1]^d; t\geq  0 \}$,  we have, for all $n$ sufficiently large,
\begin{eqnarray}\label{app}
\boldsymbol{\nabla}_{1;n}&\leq& \sup
_{\mathbf{u},\mathbf{v}\in[0,1]^d}\sup_{|\mathbf{u}-\mathbf{v}|\leq
h}|\mathbb{A}_n(\mathbf{v})-\mathbb{A}_n(\mathbf{u})|\left|\int
k(\mathbf{v})d\mathbf{v}\right|\nonumber\\&=&o_{\mathbb{P}}(1)O(1)=o_{\mathbb{P}}(1).
\end{eqnarray}
We will next evaluate $\boldsymbol{\nabla}_{2;n}$ in the right side of (\ref{decom}). We have
$$\sup_{{\mathbf{u}}\in[0,1]^d}|\mathbb{A}_n(\mathbf{u})|=O_{\mathbb{P}}(1)$$
and as $n$ tends to infinity,  by condition (C.2), we conclude that 
$$\sqrt{n}\left|\int_{\prod_{i=1}^d\left[\frac{u_i-1}{h^{1/d}},\frac{u_i}{h^{1/d}}\right]}
k(\mathbf{v})d\mathbf{v}-1\right|=o(1).$$ Then we obtain
\begin{equation}\label{DE2}
\boldsymbol{\nabla}_{2;n}=o_{\mathbb{P}}(1).
\end{equation}
As we  closely follow the lines of the proof of $\boldsymbol{\nabla}_{2;n}$, thus we obtain
\begin{equation}\label{DE4}\boldsymbol{\nabla}_{4;n}=o_{\mathbb{P}}(1).
\end{equation}
Therefore from (\ref{decom}), (\ref{bia}), (\ref{app}), (\ref{DE2}) and (\ref{DE4}), we conclude that

\begin{equation}\label{dev}
\sup_{{\bf {u}}\in[0,1]^d}| \widehat{\mathbb{A}}_n({\bf{u}})-
\mathbb{A}_n({\bf{u}})|=o_{\mathbb{P}}(1).
\end{equation}
An application of the triangle inequality shows, in turn, that
\begin{eqnarray*}
\sup_{{\bf {u}}\in[0,1]^d}\left|\widehat{\mathbb{A}}_n({\bf
u})-\frac{1}{\sqrt{n}}\mathscr{K}^*_{\mathbb{ C}}({\bf u},n)\right|&\leq&\sup_{{\bf
{u}}\in[0,1]^d}|\widehat{\mathbb{A}}_n({\bf
u})-\mathbb{A}_n({\bf u})|\\&&+\sup_{{\bf
{u}}\in[0,1]^d}\left|\mathbb{A}_n({\bf
u})-\frac{1}{\sqrt{n}}\mathscr{K}^*_{\mathbb{ C}}({\bf u},n)\right|.
\end{eqnarray*}
This, when combined with (\ref{dev}) and Theorem \ref{th}, completes the proof of Corollary  \ref{th1}.
\hfill$\Box$ \\
%\section*{Acknowledgements}
%The authors would like to thank an Associate-Editor and two referees for their very helpful comments,
%which led to a considerable improvement of the original version of the paper.


\begin{thebibliography}{}
\bibitem[Adler(1990)]{Adler1990}
Adler, R.~J. (1990).
\newblock {\em An introduction to continuity, extrema, and related topics for
 general {G}aussian processes}.
\newblock Institute of Mathematical Statistics Lecture Notes---Monograph
 Series, 12. Institute of Mathematical Statistics, Hayward, CA.


\bibitem[Bouzebda {\em et~al.}(2011a)]{Bouzebda2010}
Bouzebda, S., El~Faouzi, N.-E., and Zari, T. (2011a).
\newblock On the multivariate two-sample problem using strong approximations of
  empirical copula processes.
\newblock {\em Comm. Statist. Theory Methods}, {\bf 40}(8),
1490--1509.

\bibitem[Bouzebda {\em et~al.}(2011b)]{Bouzebda2011MMS}
Bouzebda, S., Keziou, A., and Zari, T. (2011b).
\newblock $K$-sample problem using strong approximations of empirical copula
  processes.
\newblock {\em Math. Methods Statist.}, {\bf 20}(2), 14--29.

\bibitem[B{\"u}cher and Dette(2010)]{Dette20102}
B{\"u}cher, A. and Dette, H. (2010).
\newblock A note on bootstrap approximations for the empirical copula process.
\newblock {\em Statist. Probab. Lett.}, {\bf 80} (23--24),  1925--1932.

\bibitem[Chen and Huang(2007)]{Chen200777}
Chen, S.~X. and Huang, T.-M. (2007).
\newblock Nonparametric estimation of copula functions for dependence
  modelling.
\newblock {\em Canad. J. Statist.}, {\bf 35}(2), 265--282.

\bibitem[Cherubini {\em et~al.}(2004)]{Cherubini2004}
Cherubini, U., Luciano, E., and Vecchiato, W. (2004).
\newblock {\em Copula methods in finance}.
\newblock Wiley Finance Series. John Wiley \& Sons Ltd., Chichester.

\bibitem[Chung(1949)]{sChung1949}
Chung, K.-L. (1949).
\newblock An estimate concerning the {K}olmogoroff limit distribution.
\newblock {\em Trans. Amer. Math. Soc.}, {\bf 67}, 36--50.

\bibitem[Cs{\"o}rg{\H{o}}(1979)]{csorgo1979}
Cs{\"o}rg{\H{o}}, M. (1979).
\newblock Strong approximations of the {H}oeffding, {B}lum, {K}iefer,
  {R}osenblatt multivariate empirical process.
\newblock {\em J. Multivariate Anal.}, {\bf 9}(1), 84--100.

\bibitem[Cs{\"o}rg{\H{o}} and Horv{\'a}th(1988)]{Csorgo1988}
Cs{\"o}rg{\H{o}}, M. and Horv{\'a}th, L. (1988).
\newblock A note on strong approximations of multivariate empirical processes.
\newblock {\em Stochastic Process. Appl.}, {\bf 28}(1), 101--109.

\bibitem[Cs{\"o}rg{\H{o}} and Horv{\'a}th(1993)]{Csorgoho1993}
Cs{\"o}rg{\H{o}}, M. and Horv{\'a}th, L. (1993).
\newblock {\em Weighted approximations in probability and statistics}.
\newblock Wiley Series in Probability and Mathematical Statistics: Probability
  and Mathematical Statistics. John Wiley \& Sons Ltd., Chichester.
\newblock With a foreword by David Kendall.

\bibitem[Cs{\"o}rg{\H{o}} {\em et~al.}(1997)]{Cosorgo1997}
Cs{\"o}rg{\H{o}}, M., Horv{\'a}th, L., and Szyszkowicz, B. (1997).
\newblock Integral tests for suprema of {K}iefer processes with application.
\newblock {\em Statist. Decisions}, {\bf 15}(4), 365--377.


\bibitem[Cs{\"o}rg{\H{o}} and R{\'e}v{\'e}sz(1981)]{Csorgo1981REVESZ}
Cs{\"o}rg{\H{o}}, M. and R{\'e}v{\'e}sz, P. (1981).
\newblock {\em Strong approximations in probability and statistics}.
\newblock Probability and Mathematical Statistics. Academic Press Inc.
  [Harcourt Brace Jovanovich Publishers], New York.

\bibitem[Cs{\"o}rg{\H{o}} and Hall(1984)]{Hall1983S}
Cs{\"o}rg{\H{o}}, S. and Hall, P. (1984).
\newblock The {K}oml\'os-{M}ajor-{T}usn\'ady approximations and their
  applications.
\newblock {\em Austral. J. Statist.}, {\bf 26}(2), 189--218.


\bibitem[Cui and Sun(2004)]{Cui2004}
Cui, S. and Sun, Y. (2004).
\newblock Checking for the gamma frailty distribution under the marginal
  proportional hazards frailty model.
\newblock {\em Statist. Sinica}, {\bf 14}(1), 249--267.

\bibitem[DasGupta(2008)]{Dasgupta2008}
DasGupta, A. (2008).
\newblock {\em {Asymptotic theory of statistics and probability.}}
\newblock {Springer Texts in Statistics. New York, NY: Springer.}


\bibitem[Deheuvels(1979)]{deheuvels1979a}
Deheuvels, P. (1979).
\newblock La fonction de d\'ependance empirique et ses propri\'et\'es. {U}n
  test non param\'etrique d'ind\'ependance.
\newblock {\em Acad. Roy. Belg. Bull. Cl. Sci. (5)}, {\bf 65}(6), 274--292.

\bibitem[Deheuvels(1980)]{Deheuvels1981b}
Deheuvels, P. (1980).
\newblock Nonparametric test of independence.
\newblock In {\em Nonparametric asymptotic statistics (Proc. Conf., Rouen,
  1979) (French)},  {\em Lecture Notes in Math.}, Vol. 821, J.P. Raoult, ed.,  pages 95--107.
  Springer, Berlin.

\bibitem[Deheuvels(1981)]{deheuvels1981}
Deheuvels, P. (1981).
\newblock Multivariate tests of independence.
\newblock In {\em Analytical methods in probability theory (Oberwolfach,
  1980)}, volume 861 of {\em Lecture Notes in Math.}, pages 42--50. Springer,
  Berlin.

\bibitem[Deheuvels(2009)]{Dehe100000}
Deheuvels, P. (2009).
\newblock A multivariate {B}ahadur-{K}iefer representation for the empirical
  copula process.
\newblock {\em Zap. Nauchn. Sem. S.-Peterburg. Otdel. Mat. Inst. Steklov.
  (POMI)}, {\bf 364}(Veroyatnost i Statistika. 14.2), 120--147, 237.

\bibitem[Fermanian {\em et~al.}(2004)]{fermanianradulovicdragan2004}
Fermanian, J.-D., Radulovi{\'c}, D., and Wegkamp, M. (2004).
\newblock Weak convergence of empirical copula processes.
\newblock {\em Bernoulli}, {\bf 10}(5), 847--860.

\bibitem[Frees and Valdez(1998)]{frees1998}
Frees, E.~W. and Valdez, E.~A. (1998).
\newblock Understanding relationships using copulas.
\newblock {\em N. Am. Actuar. J.}, {\bf 2}(1), 1--25.

\bibitem[Koml{\'o}s {\em et~al.}(1975)]{KMT}
Koml{\'o}s, J., Major, P., and Tusn{\'a}dy, G. (1975).
\newblock An approximation of partial sums of independent {${\rm RV}$}'s and
  the sample {${\rm DF}$}. {I}.
\newblock {\em Z. Wahrscheinlichkeitstheorie und Verw. Gebiete}, {\bf 32},
  111--131.

\bibitem[Gaenssler and Stute(1987)]{stute1987}
Gaenssler, P. and Stute, W. (1987).
\newblock {\em Seminar on empirical processes}, volume~9 of {\em DMV Seminar}.
\newblock Birkh\"auser Verlag, Basel.

\bibitem[Joe(1997)]{Joe1997}
Joe, H. (1997).
\newblock {\em Multivariate models and dependence concepts}, volume~73 of {\em
  Monographs on Statistics and Applied Probability}.
\newblock Chapman \& Hall, London.

\bibitem[McNeil {\em et~al.}(2005)]{McNeil2005}
McNeil, A.~J., Frey, R., and Embrechts, P. (2005).
\newblock {\em Quantitative risk management}.
\newblock Princeton Series in Finance. Princeton University Press, Princeton,
  NJ.
\newblock Concepts, techniques and tools.

\bibitem[Moore and Spruill(1975)]{Moore1975}
Moore, D.~S. and Spruill, M.~C. (1975).
\newblock Unified large-sample theory of general chi-squared statistics for
  tests of fit.
\newblock {\em Ann. Statist.}, {\bf 3}, 599--616.


\bibitem[Nelsen(2006)]{nelsen2006}
Nelsen, R.~B. (2006).
\newblock {\em An introduction to copulas}.
\newblock Springer Series in Statistics. Springer, New York, second edition.

\bibitem[Omelka {\em et~al.}(2009)]{Omelka2009}
Omelka, M., Gijbels, I., and Veraverbeke, N. (2009).
\newblock Improved kernel estimation of copulas: weak convergence and
  goodness-of-fit testing.
\newblock {\em Ann. Statist.}, {\bf 37}(5B), 3023--3058.


\bibitem[Philipp and Pinzur(1980)]{Philipp1980}
Philipp, W. and Pinzur, L. (1980).
\newblock Almost sure approximation theorems for the multivariate empirical
  process.
\newblock {\em Z. Wahrsch. Verw. Gebiete}, {\bf 54}(1), 1--13.

\bibitem[Piterbarg(1996)]{Piterbarg1996}
Piterbarg, V.~I. (1996).
\newblock {\em Asymptotic methods in the theory of {G}aussian processes and
 fields}, volume 148 of {\em Translations of Mathematical Monographs}.
\newblock American Mathematical Society, Providence, RI.
\newblock Translated from the Russian by V. V. Piterbarg, Revised by the
 author.
\bibitem[R{\'e}millard and Scaillet(2009)]{Remillard2009}
R{\'e}millard, B. and Scaillet, O. (2009).
\newblock Testing for equality between two copulas.
\newblock {\em J. Multivariate Anal.}, {\bf 100}(3), 377--386.

\bibitem[R{\'e}millard(2010)]{Remillard2010}
R{\'e}millard, B.
\newblock Goodness-of-Fit Tests for Copulas of Multivariate Time Series (December 22, 2010).
\newblock  {\em Available at SSRN: http://ssrn.com/abstract=1729982}

\bibitem[R{\"u}schendorf(1974)]{Ruch1}
R{\"u}schendorf, L. (1974).
\newblock On the empirical process of multivariate, dependent random variables.
\newblock {\em J. Multivariate Anal.}, {\bf 4}, 469--478.

\bibitem[R{\"u}schendorf(1976)]{Rush2}
R{\"u}schendorf, L. (1976).
\newblock Asymptotic distributions of multivariate rank order statistics.
\newblock {\em Ann. Statist.}, {\bf 4}(5), 912--923.

\bibitem[R{\"u}schendorf(2009)]{ruchendorf2009}
R{\"u}schendorf, L. (2009).
\newblock On the distributional transform, {S}klar's theorem, and the empirical
  copula process.
\newblock {\em J. Statist. Plann. Inference}, {\bf 139}(11), 3921--3927.

\bibitem[Ruymgaart(1973)]{Ruymgaart1973}
Ruymgaart, F. (1973).
\newblock {\em Asymptotic Theory for Rank Tests for Independence, MC Tract 43}.
\newblock Ph.D. thesis, Amsterdam: Mathematisch Institut.

\bibitem[Ruymgaart {\em et~al.}(1972)]{Ruymgaart_Shorack_Zwet1972}
Ruymgaart, F.~H., Shorack, G.~R., and van Zwet, W.~R. (1972).
\newblock Asymptotic normality of nonparametric tests for independence.
\newblock {\em Ann. Math. Statist.}, {\bf 43}, 1122--1135.

\bibitem[Scaillet(2005)]{Scaillet2005}
Scaillet, O. (2005).
\newblock A {K}olmogorov-{S}mirnov type test for positive quadrant dependence.
\newblock {\em Canad. J. Statist.}, {\bf 33}(3), 415--427.


\bibitem[Segers(2010)]{Segers2010}
Segers, J. (2010).
\newblock Weak convergence of empirical copula processes under nonrestrictive
  smoothness assumptions.
\newblock {\em ArXiv e-prints}.

\bibitem[Shorack and Wellner(1986)]{ShorackGalen1986s}
Shorack, G.~R. and Wellner, J.~A. (1986).
\newblock {\em Empirical processes with applications to statistics}.
\newblock Wiley Series in Probability and Mathematical Statistics: Probability
  and Mathematical Statistics. John Wiley \& Sons Inc., New York.

\bibitem[Sklar(1959)]{Sklar1959}
Sklar, A. (1959).
\newblock Fonctions de r\'epartition \`a {$n$} dimensions et leurs marges.
\newblock {\em Publ. Inst. Statist. Univ. Paris}, {\bf 8}, 229--231.

\bibitem[Sklar(1973)]{Sklar1973}
Sklar, A. (1973).
\newblock Random variables, joint distribution functions, and copulas.
\newblock {\em Kybernetika (Prague)}, {\bf 9}, 449--460.

\bibitem[Stute(1982)]{Stute1982}
Stute, W. (1982).
\newblock {The oscillation behavior of empirical processes.}
\newblock {\em Ann. Probab.}, {\bf 10}, 86--107.

\bibitem[Stute(1984)]{stute1984}
Stute, W. (1984).
\newblock {The oscillation behavior of empirical processes: The multivariate
  case.}
\newblock {\em Ann. Probab.}, {\bf 12}, 361--379.

\bibitem[Schweizer(1991)]{Schweizer1991}
Schweizer, B. (1991).
\newblock Thirty years of copulas.
\newblock In {\em Advances in probability distributions with given marginals
  (Rome, 1990)}, volume~67 of {\em Math. Appl.}, Vol. 67, G. Dall'Aglio, S. Kotz, and G. Salinetti, eds., pages 13--50. Kluwer Acad.
  Publ., Dordrecht.

\bibitem[Tsukahara(2000)]{Tsukahara2000}
Tsukahara, H. (2000).
\newblock Empirical copulas and some applications.
\newblock Research Report~27, The Institute for Economic Studies, Seijo
  University.


\bibitem[Tsukahara(2005)]{tsukahara2005}
Tsukahara, H. (2005).
\newblock Semiparametric estimation in copula models.
\newblock {\em Canad. J. Statist.}, {\bf 33}(3), 357--375.

\bibitem[van~der Vaart and Wellner(1996)]{Wellner1996}
van~der Vaart, A.~W. and Wellner, J.~A. (1996).
\newblock {\em Weak convergence and empirical processes}.
\newblock Springer Series in Statistics. Springer-Verlag, New York.
\newblock With applications to statistics.

\bibitem[Wichura(1973)]{wicura1973}
Wichura, M.~J. (1973).
\newblock Some {S}trassen-type laws of the iterated logarithm for
  multiparameter stochastic processes with independent increments.
\newblock {\em Ann. Probability}, {\bf 1}, 272--296.

\end{thebibliography}
\end{document}